\documentclass[preprint,11pt]{elsarticle}

\usepackage{graphicx}

\usepackage{amssymb, amsmath}
\usepackage{amsthm}

\usepackage{lineno}

\journal{Involve :  A journal of mathematics}

\begin{document}

\begin{frontmatter}



\title{A NEW APPROACH TO FIND AN APPROXIMATE SOLUTION OF LINEAR INITIAL VALUE PROBLEMS}


\author[label1]{Udaya Pratap Singh\corref{CorrAuthor}}
 \ead{upsingh1980@gmail.com}
 \fntext[label1]{Corresponding Author}
 \cortext[CorrAuthor]{a}
 \address[label1]{Department of Applied Sciences, Rajkiya Engineering College, Sonbhadra, Uttar Pradesh, India}


\begin{abstract}
This work investigates a new approach to find closed form analytical approximate solution of linear initial value problems. Classical Bernoulli polynomials have been used to derive a finite set of orthonormal polynomials and a finite operational matrix to simplify derivatives of dependent variable. These orthonormal polynomials together with the operational matrix of relevant order provides a good approximation to the solution of a linear initial value problem. Depending upon the nature of a problem, a series form approximation or numerical approximation can be obtained. The technique has been demonstrated through three problems. Approximate solutions have been compared with available exact or other numerical solutions. High degree of accuracy has been noted in numerical values of solutions for considered problems.
\end{abstract}

\begin{keyword}
Science \sep Publication \sep Complicated
approximate solution \sep Bernoulli polynomials \sep initial value problems \sep orthonormal polynomials

\vspace{0.7em}

\MSC 34A45 \sep B4B05 \sep 11B68
\end{keyword}

\end{frontmatter}


\section{Introduction}

Initial value problems (IVPs) for ordinary differential equations arise in a natural way in real life problems and modelling of science and engineering problems. Some examples of such modelling problems are heat conduction, wave propagation, diffusion problems, gas dynamics, nuclear physics, atomic structures, fluid flow and chemical reactions, continuum mechanics, electricity and magnetism, geophysics, antenna, synthesis problem, population genetics communication theory, mathematical modelling of economics, radiation problems and astrophysics. Many times, the exact analytic solution of such problems are not available which give rise a need to find numerical solutions. Bulk of literature is available to explore exact and numerical solutions of initial and boundary value problems \cite{Pandey2016,Atkinson1997,Samadyar2019,Bhrawy2012}. Many researchers have focused their attention to find approximate solutions differential and integral equations. Xu \cite{Xu2007} adopted method of variational iteration, Pandey, et. al. \cite{Pandey2009} applied homotopic perturbation and method of collocation. Cheon \cite{Cheon2003} discussed possible applications of Bernoulli polynomials and functions in numerical analysis. Some other latest investigations include uses of Chebyshev polynomials \cite{Maleknejad2007}, Legendre polynomials \cite{Nemati2015}, Laguerre polynomials and Wavelet Galerkin method \cite{Rahman2012}, Legendre wavelets \cite{Yousefi2006}, the operational matrix \cite{Sahu2019}. 

Numerical solution to boundary value problems has been also presented by Shiralashetti and Kumbinarasaiah \cite{Shiralashetti2018,Shiralashetti2019}, 	Abd-Elhameed et al. \cite{Abd-Elhameed2013}, Iqbal et al. \cite{Iqbal2018} and, Kumar and Singh \cite{KumarM2009}. Bernoulli polynomials and its properties have been also discussed by many authors \cite{Kurt2011,Natalini2003}. Tohidi et. al. \cite{Tohidi2013} obtained numerical approximation for generalized pantograph equation using Bernoulli matrix method, Tohidi and Khorsand \cite{Tohidi2013a} to solve second-order linear system of partial differential equations, Mohsenyzadeh \cite{Mohsenyzadeh2016} used Bernoulli polynomials to solve Volterra type integral equations. Recently, Singh et al. \cite{Singh2019} used Bernoulli polynomials to develop a trigonal operational matrix to solve Abel-Volterra type integral equations. 

In this work, it is proposed to solve linear initial value problems using orthogonal polynomials derived from Bernoulli polynomials with a modified operational matrix \cite{Singh2019}.

\section{Bernoulli Polynomials}
The word Bernoulli Polynomials was first coined by J. L. Raabe in 1851 while discussing the formula $\sum\limits_{n = 0}^{m - 1} {{B_n}\left( {x + \frac{k}{m}} \right)}  = {m^{ - (n + 1)}}{B_n}(mx)$, however, the polynomials $B_n(x)$ were already introduced by Jakob Bernoulli in 1690 in his book "$Ars Conjectandi$" \cite{Costabile2006}. A thorough study of these polynomials was first done by Leonhard Euler in 1755, who showed in his book “Foundations of differential calculus” that these polynomials satisfy the finite difference relation:
\begin{equation} \label{eq. 1 : rec. rel. of Bn} 
{B_n}(\zeta  + 1) - {B_n}(\zeta ) = n{\zeta ^{n - 1}},\,\,\,\,\,n \ge 1
\end{equation}
and proposed the method of generating function to calculate  $B_n(x)$. Following Leonhard Euler, recently Costabile and Dell’Accio \cite{Costabile2006} showed that Bernoulli Polynomials are monic which can be extracted from its generating function  
\begin{equation} \label{eq.2 : gen. fn. of Bn}
\frac{{\gamma \,{e^{\zeta \gamma }}}}{{{e^\gamma } - 1}} = \sum\limits_{n = 0}^\infty  {{B_n}(\zeta )\,\frac{{{\gamma ^n}}}{{n!}}} \,\,\,\,\,\left( {\left| \zeta  \right| < 2\pi } \right)
\end{equation}
and represented in the simple form:
\begin{equation}\label{eq.1 : Basic Bernoulli Polynomials}
{B_n}(\zeta ) = \sum\limits_{j = 0}^n {\,\left( {\begin{array}{*{20}{c}}
		n \\ j 	\end{array}} \right)\,\,{B_j}(0)\,{\zeta ^{n - j}},\hspace{1em} n = 0,1,2,...} \hspace{1em} 0 \le \zeta  \le 1
\end{equation}
where, $B_n(0)$ are the Bernoulli numbers, which can also be calculated with Kronecker’s formula ${B_n}(0) =  - \sum\limits_{j = 1}^{n + 1} {\frac{{{{( - 1)}^j}}}{j}\left( {\begin{array}{*{20}{c}}
{n + 1}\\j \end{array}} \right)} \sum\limits_{k = 1}^j {{k^n}} \,\,;\,\,\,n \ge 0$ \cite{Todorov1984}. Thus, first few Bernoulli polynomials can be written as ${B_0}(\zeta ) = 1, B_1(\zeta)=\zeta-\frac{1}{2}, B_2(\zeta)=\zeta^2-\zeta+\frac{1}{6}	, 	B_3(\zeta)=\zeta^3-\frac{3}{2}\zeta^2+\frac{1}{2}\zeta	, 	B_4(\zeta)=\zeta^4-2\zeta^3+\zeta^2-\frac{1}{30}$.

These Bernoulli Polynomials form a complete basis over $[0,1]$ \cite{Kreyszig1978} and show some interesting properties \cite{Costabile2001,Lu2011}:
\begin{equation}\label{eq.4 : properties of Bernoulli Polnmls}
\left. \begin{gathered}
B'_n(\zeta ) = n{B_{n - 1}}(\zeta ),\,\,\,n \ge 1 \hfill \\
\int_0^1 {{B_n}(z)dz = 0,\,\,\,\,\,\,\,\,\,\,n \ge 1}  \hfill \\
{B_n}(\zeta  + 1) - {B_n}(\zeta ) = n{\zeta ^{n - 1}},\,\,n \ge 1 \hfill \\ 
\end{gathered}  \right\}.
\end{equation}

\section{The Orthonormal Polynomials}
It can be easily verified that the polynomials ${B_n}(x)\,(n \ge 1)$  given by (eq.\ref{eq.1 : Basic Bernoulli Polynomials})  are orthogonal to $B_o(x)$ with respect to standard inner product on  ${L^2} \in [0,1]$. Using this property, an orthonormal set of polynomials can be derived for any $n$ with Gram-Schmidt orthogonalization. First ten such orthonormal polynomials are obtained  in \ref{Appendix A}.

\section{Approximation of Functions}

\hspace{-2em} \textbf{Theorem.} Let $H=L^2[0,1]$  be a Hilbert space and $Y=span\left\{y_0,y_1,y_2,...,y_n\right\}$  be a subspace of $H$ such that $dim{(}Y)<\infty$ , every $f\in H$  has a unique best approximation out of $Y$ \cite{Kreyszig1978}, that is, $\forall y(t)\in Y,\, \exists \, \hat{f}(t)\in Y$ s.t. $\parallel f(t)-\hat{f}(t)\parallel_2\le\parallel f(t)-y(t)\parallel_2$. This implies that,  $\forall \, y(t)\in Y, <f(t)-f(t), y(t)>= 0$, where $<,>$  is standard inner product on $L^2\in[0,1]$ (\textit{c.f. Theorems 6.1-1 and 6.2-5, Chapter 6} \cite{Kreyszig1978}).

\hspace{-2em} \textbf{Remark.} Let $Y=span\left\{\phi_0,\phi_1,\phi_2,...,\phi_n\right\},$ where $\phi_k\in L^2[0,1]$ are orthonormal Bernoulli polynomials. Then, from the above theorem, for any function $ f\in L^2[0,1],$
\begin{equation} \label{eq.5 : approx theorem}
f\approx\hat{f}=\sum_{k=0}^{n}{c_k\phi_k}, 
\end{equation}
where $c_k=\left\langle f,\phi_k\right\rangle,$ and $<,>$ is the standard inner product on $L^2\in[0,1]$.
For numerical approximation, series (5) can be written as:
\begin{equation} \label{eq.6 : approximation series}
f(\zeta)\simeq\sum_{k=0}^{n}{c_k\phi_k=C^T\phi(\zeta)} 
\end{equation}
where $C=\left(c_0,c_1,c_2,...,c_n\right), \phi(\zeta)=\left(\phi_0,\phi_1,\phi_2,...,\phi_n\right)$ are column vectors, and number of polynomials $n$ can be chosen to meet required accuracy.

\section{Construction of operational matrix}
The orthonormal polynomials, as derived in \ref{Appendix A}, can be expressed as:
\begin{equation}\label{eq.7 : Ortho Pol. of deg. 0}
\int_0^\zeta  {{\phi _o}(\eta )d\eta }  = {\phi _o}(\zeta ) + \frac{1}{{2\sqrt 3 }}{\phi _1}(\zeta )
\end{equation}

\begin{equation}\label{eq.8 : Ortho Pols.}
\begin{array}{l}
\int\limits_0^\zeta  {{\phi _i}(x)dx = } \,\,\,\,\,\frac{1}{{2\sqrt {(2i - 1)(2i + 1)} }}{\phi _{i - 1}}(\zeta )\\
\hspace{2cm} + \frac{1}{{2\sqrt {(2i + 1)(2i + 3)} }}{\phi _{i + 1}}(\zeta ),\,\,\,\,(\,for{\rm{ }}\,i = 1\,,2,...\,,n)
\end{array}
\end{equation}

Relations (\ref{eq.7 : Ortho Pol. of deg. 0}-\ref{eq.8 : Ortho Pols.}) can be represented in closed form as:
\begin{equation}\label{eq.9 : combined expression for Ortho Pols.}
\int\limits_0^\zeta  {\phi(\eta )d\eta  = \,\,} \Theta \,\phi(\zeta )
\end{equation}

where $\zeta\in[0,1]$ and $\Theta$ is operational matrix of order $(n+1)$ given as :
\begin{equation} \label{eq.10 : operational matrix}
\Theta \, = \frac{1}{2}\left[ {\begin{array}{*{20}{c}}
	1&{\frac{1}{{\sqrt {1.3} }}}&0& \cdots &0\\
	{ - \frac{1}{{\sqrt {1.3} }}}&0&{\frac{1}{{\sqrt {3.5} }}}& \cdots &0\\
	0&{ - \frac{1}{{\sqrt {3.5} }}}&0& \ddots & \vdots \\
	\vdots & \vdots & \ddots & \ddots &{\frac{1}{{\sqrt {\left( {2n - 1} \right)\left( {2n + 1} \right)} }}}\\
	0&0& \cdots &{ - \frac{1}{{\sqrt {\left( {2n - 1} \right)\left( {2n + 1} \right)} }}}&0
	\end{array}} \right]
\end{equation}

\section{Solution of Initial Value Problems}

Consider the linear IVP:
\begin{equation}\label{eq.11 : Gen. IPV}
\frac{{{d^2}y}}{{d{\zeta ^2}}} + P(\zeta)\frac{{dy}}{{d\zeta }} + Q(\zeta)y = r(\zeta ),\,\,\,y(0) = \alpha, \left(\frac{dy}{d\zeta}\right)_{\zeta =0} = \beta
\end{equation}

where $y=y(\zeta); p(\zeta)$ and $q(\zeta)$ and $r(\zeta)$ are continuous functions defined $[0,1]$. It is further assumed that eq. \eqref{eq.11 : Gen. IPV}  admits a unique solution on $[0,1]$, otherwise, a suitable transformation $z=z(x)$ may be applied to change the domain of $y, p, q$ and $r$.

\subsection{\textbf{Case 1 : Coefficients of $y$ and $\frac{dy}{d\zeta}$ are constant}}
\label{Case 1}
\vspace{0.3em}
In this case, taking $p$ for $P(x)$ and $q$ for $Q(x)$, eq. (\ref{eq.11 : Gen. IPV}), is written as:

\begin{equation}\label{eq.11 : Gen. IPV Case 1}
\frac{{{d^2}y}}{{d{\zeta ^2}}} + p\frac{{dy}}{{d\zeta }} + q(\zeta)y = r(\zeta ),\,\,\,y(0) = \alpha, \left(\frac{dy}{d\zeta}\right)_{\zeta =0} = \beta
\end{equation}

 Let $R=\left(r_o,r_1,...,r_n \right)$ be a real column vector such that the function $r(\zeta)$ can be approximated in terms of first $n+1$ orthonormal Bernoulli polynomials as:
\begin{equation}\label{eq.12 : r for IVP}
r\left(\zeta\right)=R \, \phi(\zeta)
\end{equation}

Let $C=\left(c_o,c_1,c_2,...,c_n\right)-$ be a column vector of $n+1$ unknown quantities. Taking
\begin{equation}\label{eq.13 : y'' for IVP}
\frac{{{d^2}y}}{{d{\zeta ^2}}} = C^T\,\phi(\zeta )
\end{equation}
eq. \eqref{eq.11 : Gen. IPV} can be re-written as:
\begin{equation}\label{eq.14 for Substitution in IVP}
{C^T}{\phi}(\zeta ) + \,\,p\,\,{C^T}\,\Theta \,{\phi}(\zeta ) + \,\,q\,\,{C^T}\,{\Theta ^2}\,{\phi}(\zeta ) = R^T\,{\phi}(\zeta )
\end{equation}
which gives, 
\begin{equation}\label{eq.15 : B^T for IVP}
C^T\,\, = {\left[ {I + p\,\,\Theta + q\,\,{\Theta ^2}} \right]^{ - 1}}\,R^T
\end{equation}
 Substituting eq. (\ref{eq.15 : B^T for IVP}) back into eq. (\ref{eq.13 : y'' for IVP}), an approximation for $y(\zeta)$  can be obtained as:
\begin{equation} \label{Eq.16 : Soln. y for IVP}
y(\zeta ) = {C^T}\,\Theta^2\,{\phi}(\zeta ).
\end{equation}

\subsection{\textbf{Case 2 : Coefficients of $y$ and $\frac{dy}{d\zeta}$ are functions of independent variable}}
\label{Case 2}
\vspace{0.3em}

Taking $P(\zeta)=a^T \phi(\zeta)$ and $Q(\zeta)=b^T\phi(\zeta)$ together with eqs. (\ref{eq.12 : r for IVP}-\ref{eq.13 : y'' for IVP}), eq. (\ref{eq.11 : Gen. IPV}) can be written as:

\begin{equation}\label{Gen case Substitution in IVP}
{C^T}{\phi}(\zeta ) + \left( a^T \phi(\zeta )\right) \left(  C^T\Theta \phi(\zeta ) \right)  + \left( b^T \phi(\zeta )\right) \left( C^T\Theta ^2\phi(\zeta )\right)  = R^T\,{\phi}(\zeta )
\end{equation}

Because $a^T \phi(\zeta )$ and $b^T \phi(\zeta )$ in second and third terms respectively on left side of eq. (\ref{Gen case Substitution in IVP}) are just the polynomials or degree $2n$, eq. (\ref{Gen case Substitution in IVP}) can be re-written as, 

\begin{equation}\label{at bt manageing}
{C^T}{\phi}(\zeta ) +  C^T\Theta \left[ \phi(\zeta ) \left( a^T \phi(\zeta )\right) \right]  +  C^T\Theta^2 \left[ \phi(\zeta ) \left( b^T \phi(\zeta )\right) \right]   = R^T\,{\phi}(\zeta )
\end{equation}

Here, $\phi(\zeta ) \left( a^T \phi(\zeta ) \right) $ is a vector of type 
\[\left(  \phi_o(\zeta )\sum_{k=0}^{k=n} a_k  \phi_k(\zeta ),\,\phi_1(\zeta )\sum_{k=0}^{k=n} a_k  \phi_k(\zeta ),\,\dots,\,\phi_n(\zeta )\sum_{k=0}^{k=n} a_k  \phi_k(\zeta ) \right) \]

\begin{equation} \label{Merging of aT & phi }
 \equiv (\psi_o(\zeta),\psi_1(\zeta),\dots,\psi_n(\zeta)) = \psi(\zeta), \hspace{5mm}(say!)
\end{equation} 

In eq. (\ref{Merging of aT & phi }), each $\psi_k(\zeta)$ can be approximated as a linear combination of orthonormal polynomials in the form $ \psi_k(\zeta) = A_k^T \phi(\zeta)$, where $A_k^T$ are vectors of form $1\times (n+1)$ for $ k = 1,2,\dots,n$, and, therefore,  $\phi(\zeta ) \left( a^T \phi(\zeta ) \right) = \psi(\zeta) = A \phi(\zeta)$, where $A = \left( A_o^T,A_1^T, \dots, A_n^T \right)_{(\underline{n+1}\times1)}$. Similarly,  $\phi(\zeta ) \left( b^T \phi(\zeta ) \right)$ can be approximated as $B^T \phi(\zeta) $ for some vector $B^T = (B_o,B_1,\dots , B_n)_{(\underline{n+1}\times1)}$ such that $B_k^T$ are real vectors of form $1\times(n+1)$. With these intermediate approximations, eq. (\ref{at bt manageing}) can be written as :

\begin{equation}\label{at bt managed}
{C^T}{\phi}(\zeta ) +  C^T\Theta A \phi(\zeta )  + C^T\Theta^2 B \phi(\zeta )   = R^T\,{\phi}(\zeta )
\end{equation}

From eq. (\ref{at bt managed}), the required coefficient vector $C$ is obtained as:
\begin{equation}
C^T = R^T \left( I + \theta A + \theta^2 B \right)^{-1} ,
\end{equation}

where, $I$ is identity matrix of order $n$. The expression for $y(\zeta)$ is obtained as:
\begin{equation}
y(\zeta ) = {C^T}\,\Theta^2\,{\phi}(\zeta )
\end{equation}

\section{Numerical Examples}
In order to discuss and establish the accuracy and efficacy of the present method, following examples have been taken. 

\subsection*{Example 1:} Let us consider the IVP  

\begin{equation} \label{Ex.1}
\frac{d^2y}{dx^2}+5 \frac{dy}{dx}+3 y = e^{-x} ; \hspace{1em} y(0)=\left( \frac{dy}{dx}\right) _{x=0} =0
\end{equation}
which has exact solution $y(x) = e^{-\frac{5}{2} x} \left( \cosh(\frac{\sqrt{13}}{2} x) + \frac{3}{\sqrt{13}}\sinh(\frac{\sqrt{13}}{2} x) \right)  - e^{-x}$.

Comparing eq.(\ref{Ex.1}) with eq. (\ref{eq.11 : Gen. IPV}) and taking $m=6$, equations (\ref{eq.12 : r for IVP} - \ref{eq.15 : B^T for IVP}) yield
\begin{equation}
{B^T} = \left( \begin{array}{cccc}
0.08086, &  -0.02459,  &  -0.00156,  &  -0.00240,  \\
-0.001046, & -0.000425, & -0.00017,  & 0
\end{array}
  \right)
\end{equation}
Using value of $B^T$ in eq. (\ref{Eq.16 : Soln. y for IVP}), an approximate solution is obtained as:
\begin{equation*}
y(x) \approx - 0.00009 x + 0.25057 x^2 - 0.08655 x^3 + 0.07415 x^4
\end{equation*}
\begin{equation}
\hspace{-2em} - 0.10716 x^5 + 0.0930 x^6 - 0.03382 x^7
\end{equation}

\begin{figure}[h]
	\centering
	\includegraphics[width=0.47\linewidth]{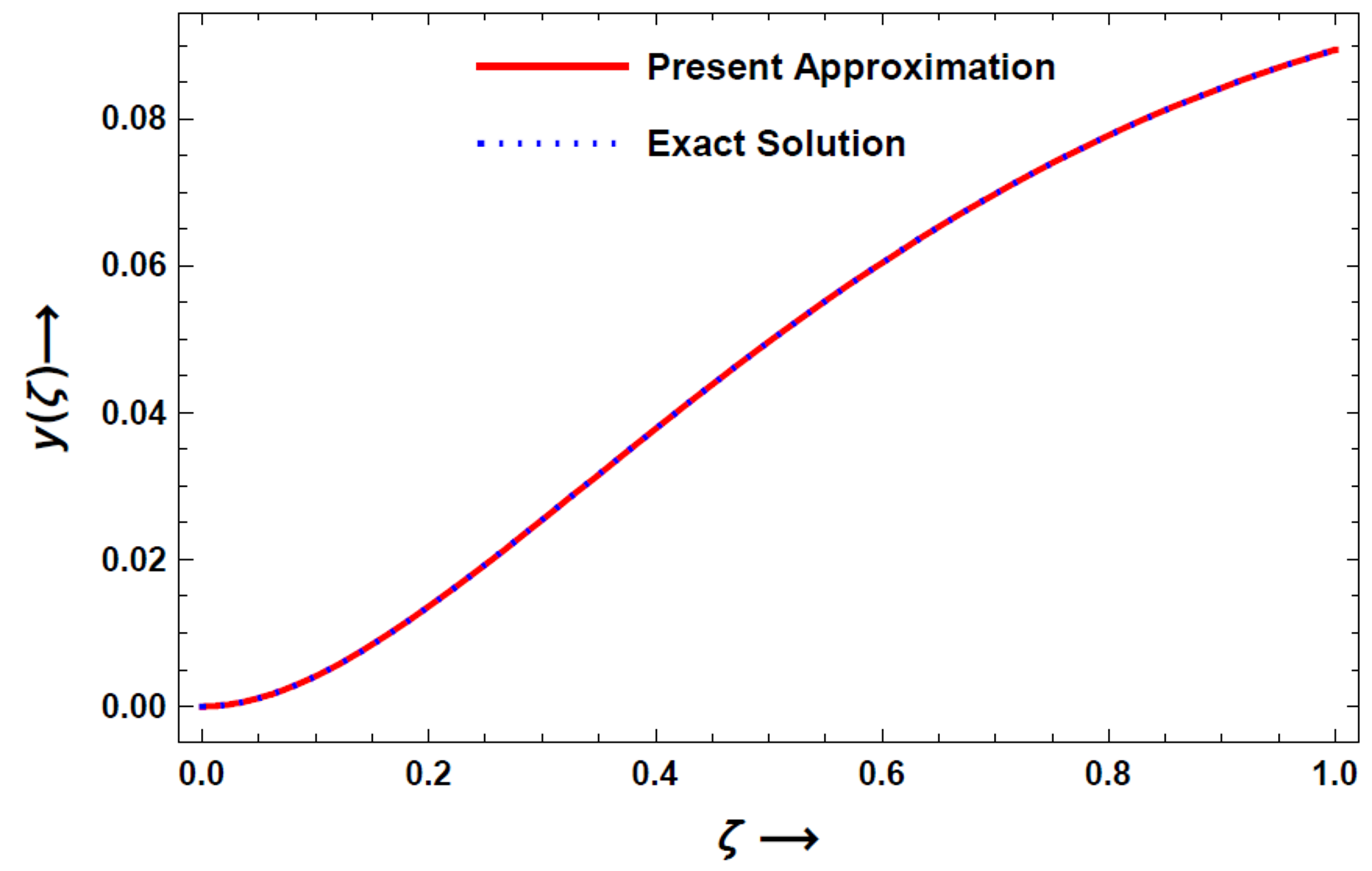} 
	\hspace{1mm}
	\includegraphics[width=0.49\linewidth]{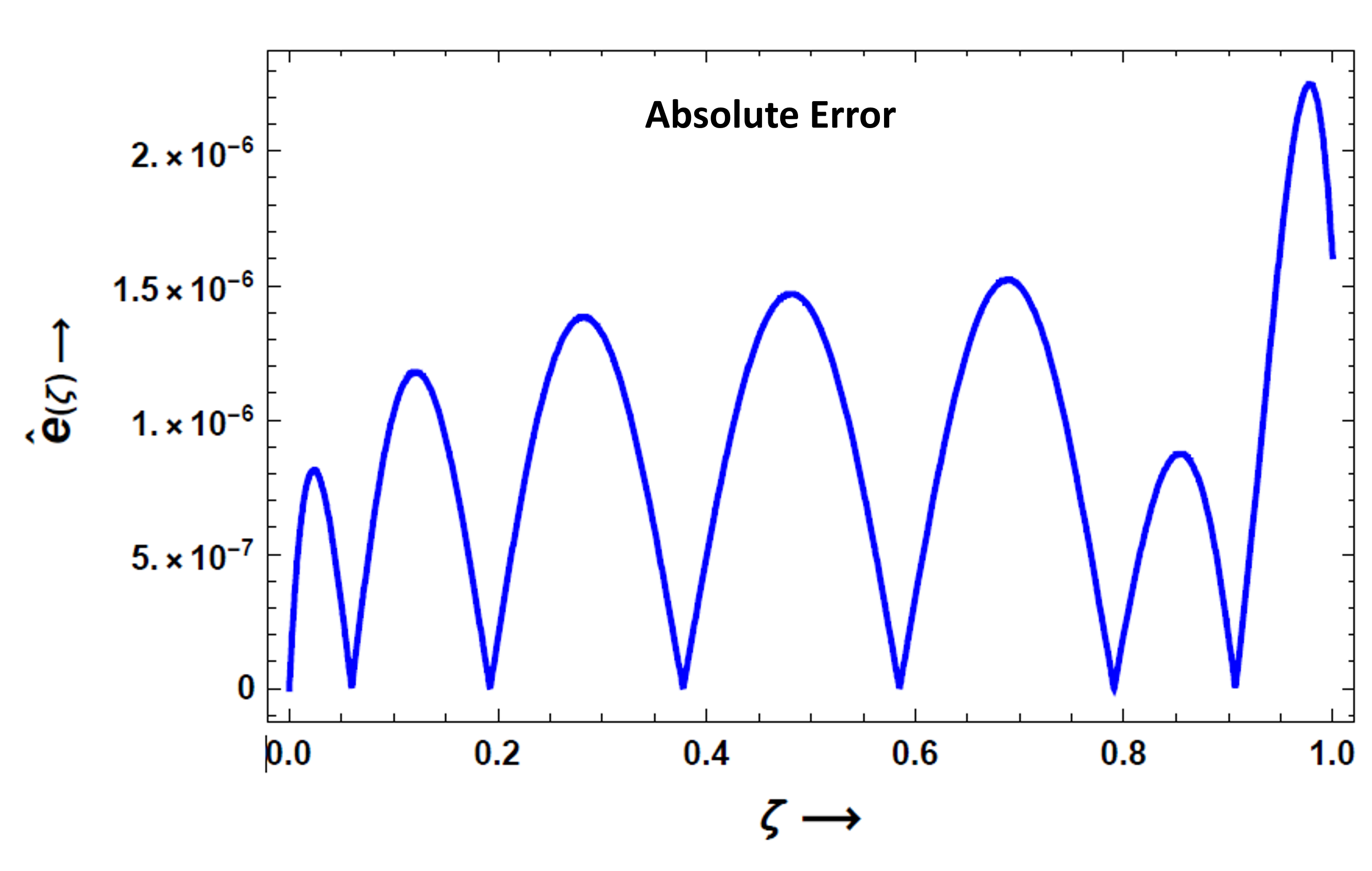}
	(a) \hspace{0.45\linewidth} (b)
	\caption{(a) Comparison of exact and present solution for example 1. (b) Absolute error between exact and approximate solutions of example 1}
	\label{fig:fig1}
\end{figure}

\subsection*{Example 2:} Consider the IVP  
\begin{equation} \label{Ex.2}
\frac{d^2y}{dx^2}-5 \frac{dy}{dx}+2 y = \tan(x) ; \hspace{1em} y(0)=\left( \frac{dy}{dx}\right) _{x=0} =0
\end{equation}
which is linear in nature but its not easy to solve manually. We will compare the present solution of this IVP with the one generated by Mathematica.

Comparing eq.(\ref{Ex.2}) with eq. (\ref{eq.11 : Gen. IPV}) and taking $m=9$, equations (\ref{eq.12 : r for IVP} - \ref{eq.15 : B^T for IVP}) yield
\begin{equation}
{B^T} = \left(  \begin{array}{c}
5.1220,\,\, 5.5181,\,\, 2.9304, \,\, 1.0668, \,\, 0.2958,\\ 0.0663,  \,\, 0.0125, \,\,  0.0020, \,\, 0.0003, \,\,  0, \,\, 0
\end{array}
\right)
\end{equation}
Using value of $B^T$ in eq. (\ref{Eq.16 : Soln. y for IVP}), an approximate solution is obtained as:
\begin{equation*}
	y(x) \approx 0.0001 x - 0.0025 x^2 + 0.1942 x^3 + 
	0.04799 x^4 + 0.7521 x^5 
\end{equation*}
\begin{equation}
\hspace{-2em} - 0.9599 x^6 + 1.5043 x^7 - 
0.9351 x^8 + 0.3669 x^9
\end{equation}

\begin{figure}[h]
	\centering
	\includegraphics[width=0.46\linewidth]{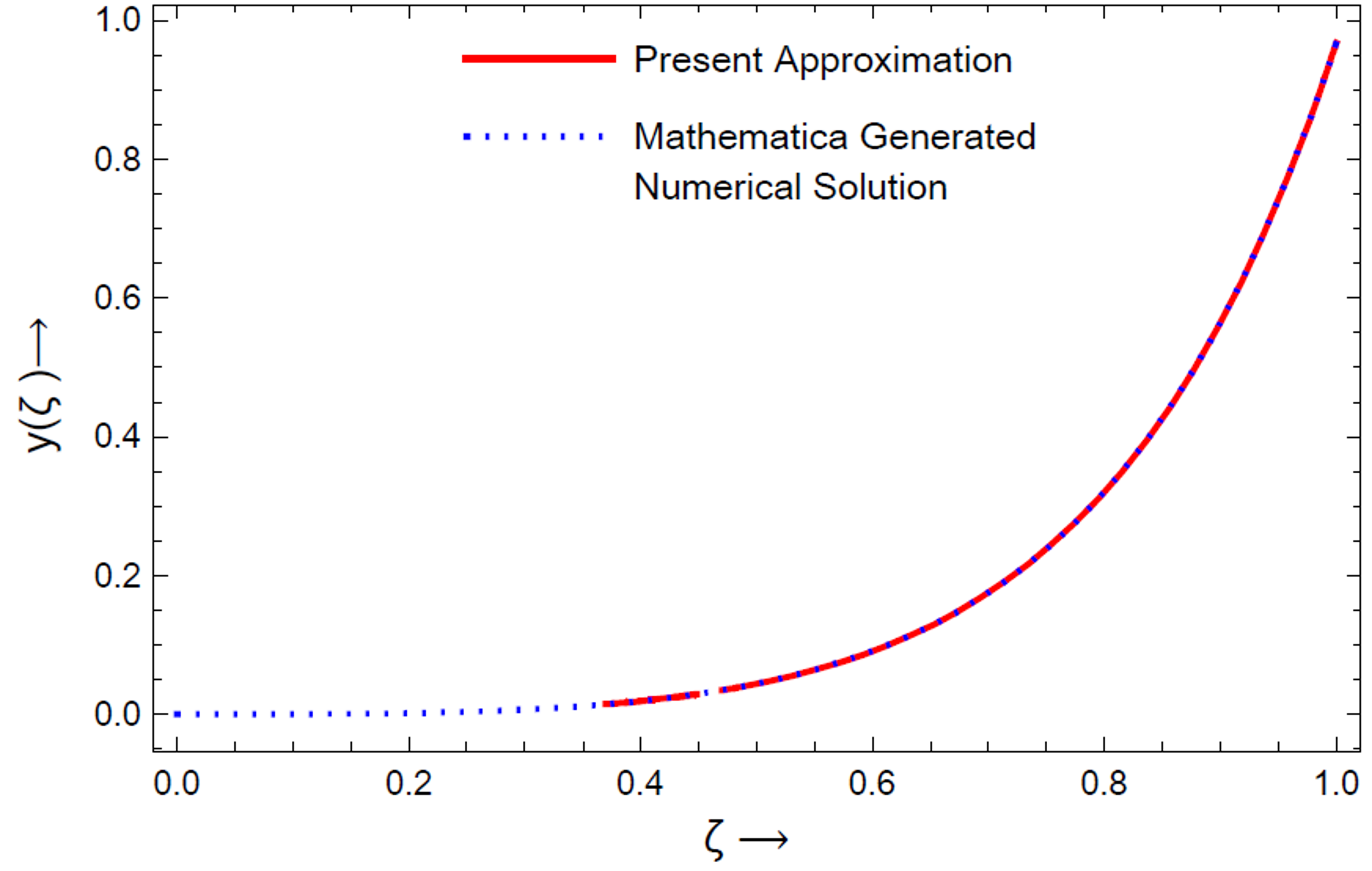} 
	\hspace{1mm}
	\includegraphics[width=0.49\linewidth]{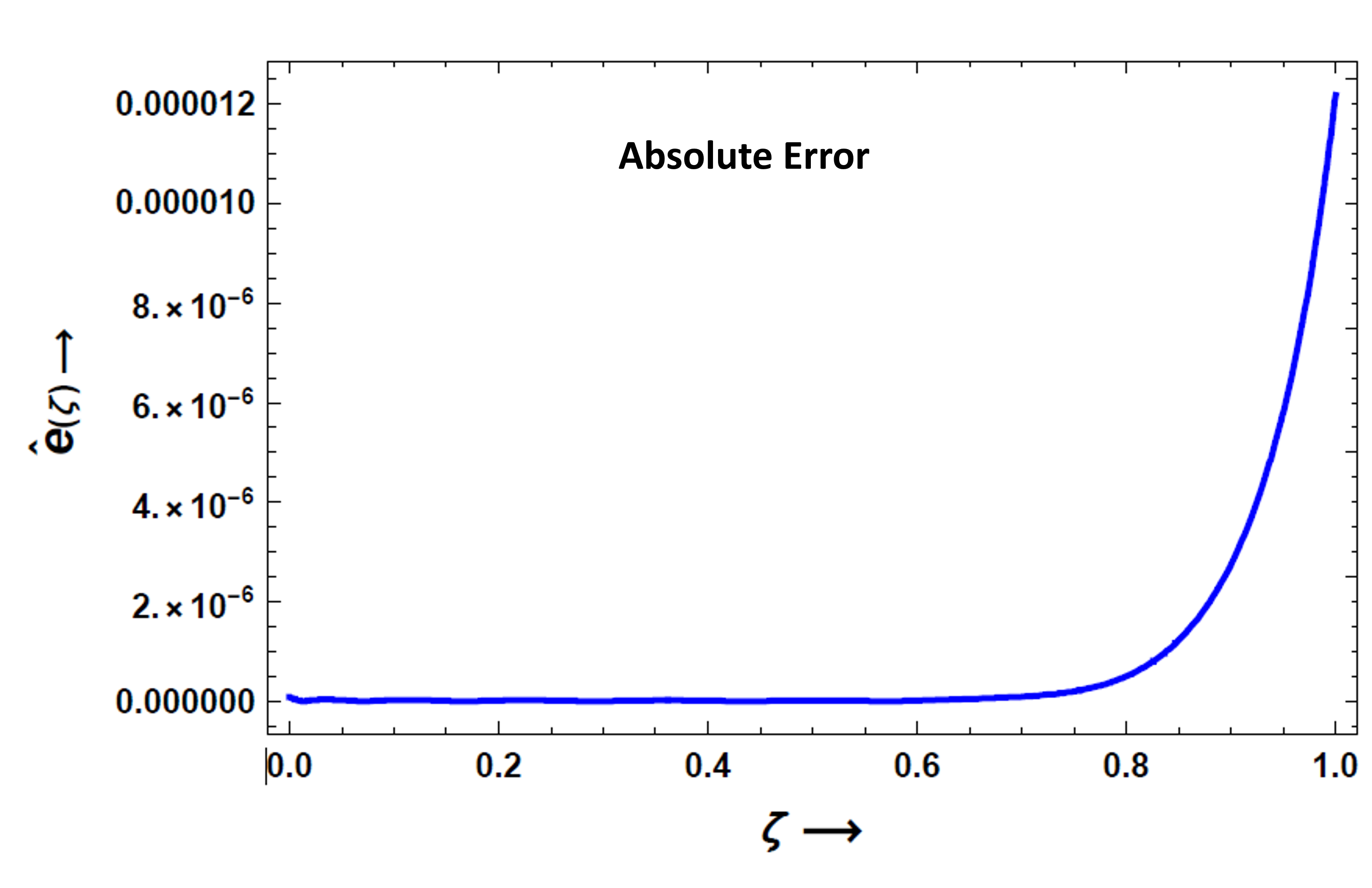}
	(a) \hspace{0.45\linewidth} (b)
	\caption{(a) Comparison of present approximation and $Mathematica$ generated numerical solutions to example 2. (b) Absolute error between present approximation and $Mathematica$ generated numerical solutions to example 2.}
	\label{fig:fig2}
\end{figure}

\subsection*{Example 3:} Let us take the IVP  
\begin{equation} \label{Ex.3}
\frac{d^2y}{dx^2}+\tan(x) \frac{dy}{dx}+2 \cos^2(x) y = 2 \cos^4(x) ; \hspace{1em} y(0)=\left( \frac{dy}{dx}\right) _{x=0} =0
\end{equation}
The exact solution of this IVP is $y(x) = 2 - 2 \cos(\sqrt{2} \sin^2 x) -\sin x $. 

Using the method discussed in $section$ \ref{Case 2}, coefficient vector $C^T$ and approximate solution $y(x)$ of example (\ref{Ex.3}) is obtained for $n=6$ as:
\begin{equation}
C^T = \left( 0.78730, 0.62821, 0.10352, -0.02660, 0.00101, 0.00136, 0.00011 \right)
\end{equation}

\begin{equation*}
	y(x) \approx -0.00025 + 0.00718 x + 2.94625 x^2 + 0.16696 x^3
\end{equation*}
\begin{equation} \label{sol. ex. 3}
\hspace{-1em} - 
1.57951 x^4 + 0.17065 x^5 + 0.33315 x^6
\end{equation}

\begin{figure}[h]
	\centering
	\includegraphics[width=0.47\linewidth]{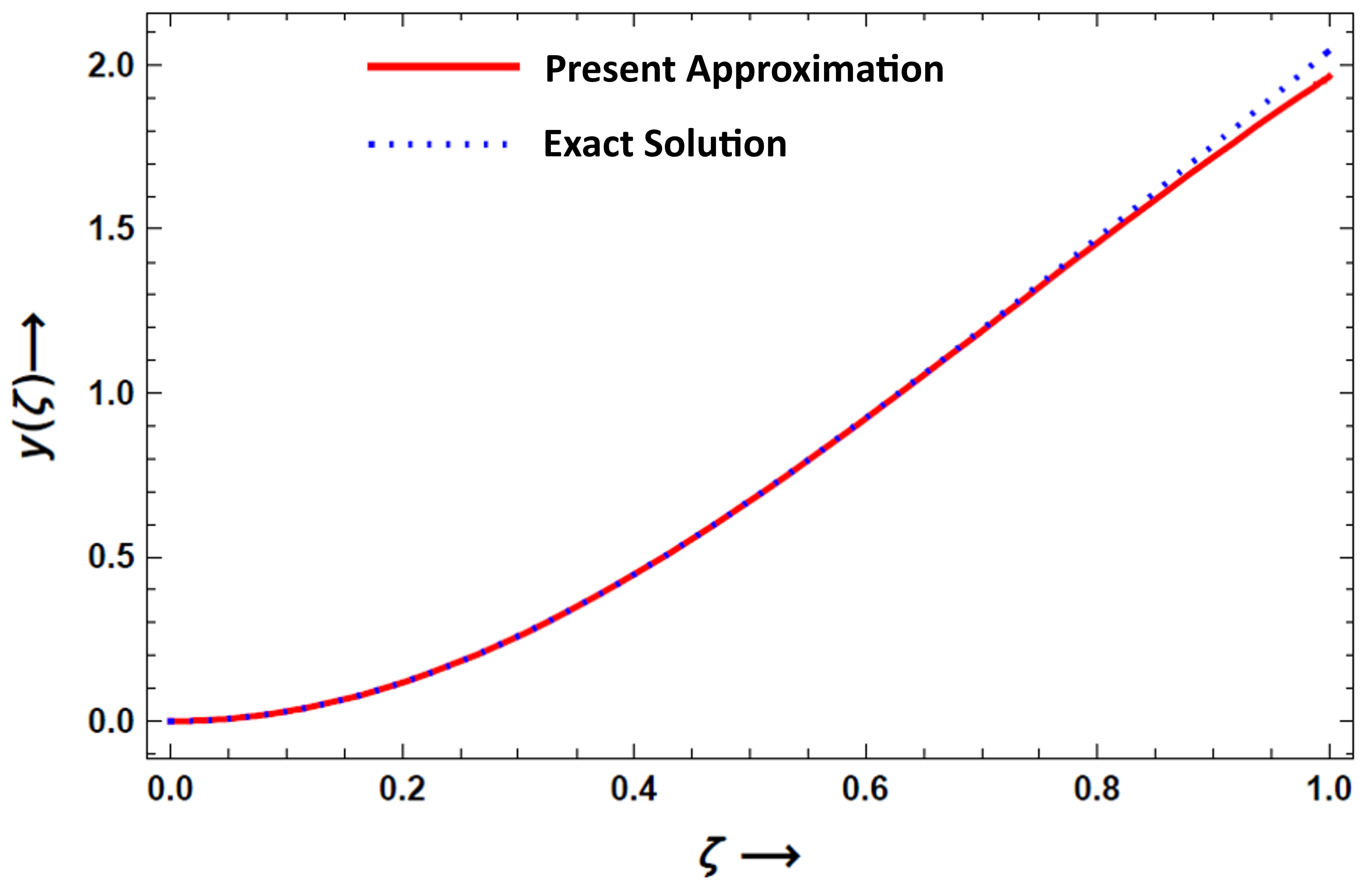} 
	\hspace{1mm}
	\includegraphics[width=0.49\linewidth]{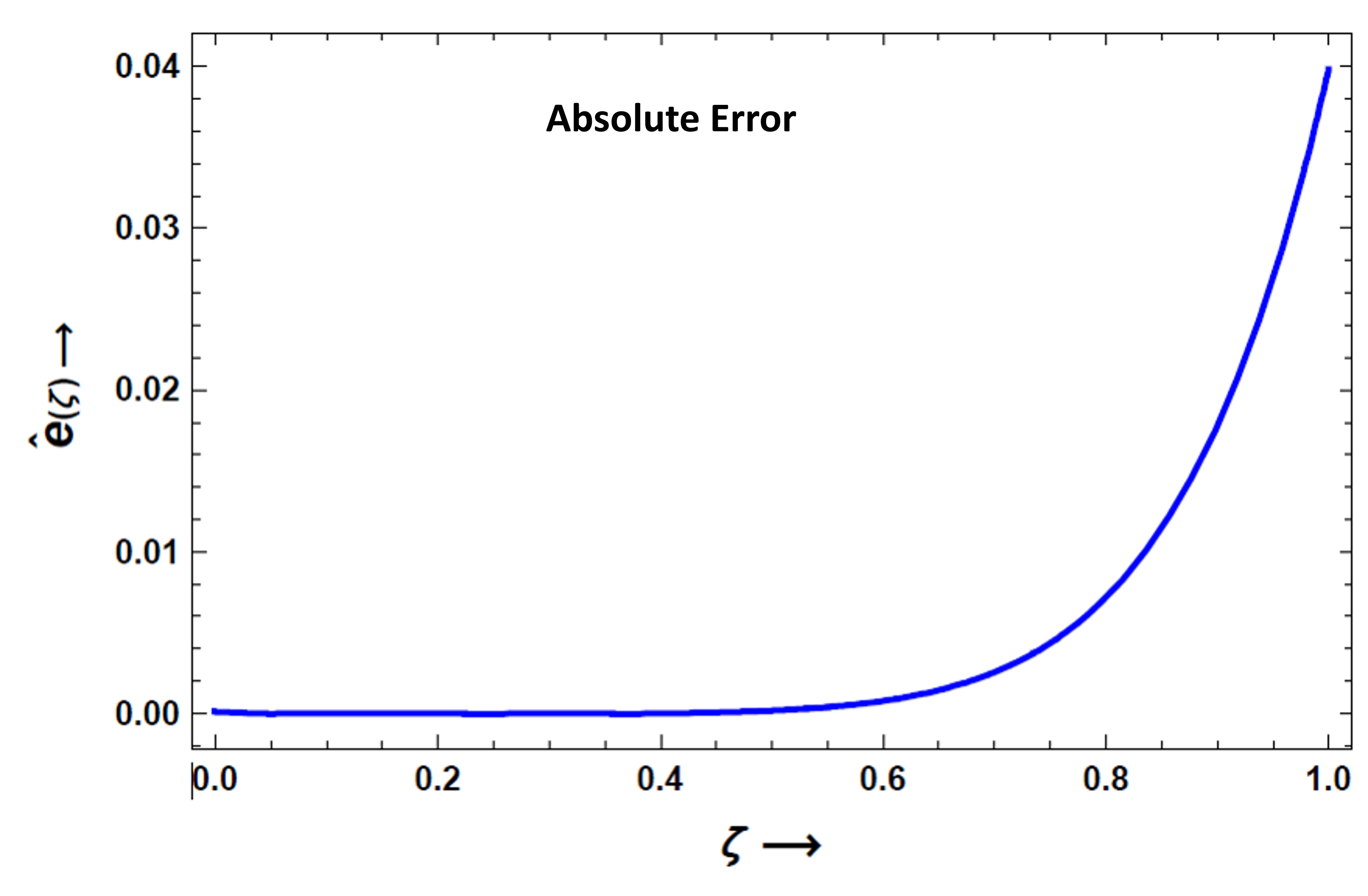}
	(a) \hspace{0.45\linewidth} (b)
	\caption{(a) Comparison of exact solution and present approximation to example 3. (b) Absolute error between exact solution and present approximation to example 3.}
	\label{fig:fig3}
\end{figure}

\section{Conclusion}
In this work, a new method was presented and demonstrated to find fast and approximate solution of linear initial value problems with help of orthogonal Bernoulli polynomials. The method includes derivation of a set of $n$ orthonormal polynomials derived from Bernoulli polynomials up to degree $n$, and an operational matrix.  The present method converts a given initial value problem into a system of algebraic equations with unknown coefficients, which are easily obtained with the help of operational matrix, and finally an approximate solution is  obtained in form of a polynomial of degree $n$. The method has been demonstrated with three examples. The main features of this method can be summarized as:
\begin{itemize}
	\item the method is programmable.
	\item solutiion is obtained in form of a polynomial of degree $n$ which can be easily used for various applications. 
	\item error can be minimized up to required accuracy because error decreases quickly with increase of $n$-the degree of Bernoulli polynomials.
	\item error is negligible for simple IVPs with constant coefficients. 
\end{itemize}




\appendix
\section{}
\label{Appendix A}

First ten orthonormal polynomials derived with Bernoulli polynomials are: 
\begin{equation}
{\phi _{0\,}}(\zeta ) = 1
\end{equation}
\begin{equation}
{\phi _1}(\zeta ) = \sqrt 3 ( - 1 + 2\zeta )
\end{equation}
\begin{equation}
\phi_2\left(x\right)=\sqrt5\left(1-6x+6x^2\right)
\end{equation}
\begin{equation}
\phi_3(\zeta)=\sqrt7(-1+12\zeta-30\zeta^2+20\zeta^3)
\end{equation}
\begin{equation}
\phi_4(\zeta)=3(1-20\zeta+90\zeta^2-140\zeta^3+70\zeta^4)
\end{equation}
\begin{equation}
\phi_5(\zeta)=\sqrt{11}(-1+30\zeta-210\zeta^2+560\zeta^3-630\zeta^4+252{\zeta5}^)
\end{equation}
\begin{equation}
{\phi _6}(\zeta ) = \sqrt {13} \left( \begin{array}{l}
1 - 42\zeta  + 420{\zeta ^2} - 1680{\zeta ^3}\\
+ 3150{\zeta ^4} - 2772\zeta {}^5 + 924{\zeta ^6}
\end{array} \right)
\end{equation}
\begin{equation}
{\phi _7}\left( x \right) = \sqrt {15} \left( \begin{array}{l}
- 1 + 56x - 756{x^2} + 4200{x^3}\\
- 11550{x^4} + 16632{x^5} - 12012{x^6} + 3432{x^7}
\end{array} \right)
\end{equation}
\begin{equation}
{\phi _8}\left( x \right) = \sqrt {17} \left( \begin{array}{l}
- 1 + 72x - 1260{x^2} + 9240{x^3} - 34650{x^4}\\
+ 72072{x^5} - 84084{x^6} + 51480{x^7} - 12870{x^8}
\end{array} \right)
\end{equation}
\begin{equation}
{\phi _9}\left( x \right) = \sqrt {19} \left( \begin{array}{l}
- 1 + 90x - 1980{x^2} + 18480{x^3} - 90090{x^4} + 252252{x^5}\\ - 420420{x^6} + 411840{x^7} - 218790{x^8} + 48620{x^9}
\end{array} \right)
\end{equation}

\bibliographystyle{model1-num-names}
\bibliography{references}







\end{document}